%% file: sixvertexfinal.tex
\documentclass[12pt,a4paper]{amsart}
\usepackage{amssymb,stmaryrd,blkarray}
\usepackage[a4paper,left=2.25cm,right=2.25cm,top=3cm,bottom=3cm,headsep=1cm]{geometry}
\usepackage{tikz}\usetikzlibrary{graphs,quotes,fit,positioning,matrix,calc,decorations.markings,angles,decorations.pathmorphing,decorations.pathreplacing}
\usepackage{hyperref}
\title[Eulerian orientations and the six-vertex model on
    planar maps]{Eulerian orientations and the six-vertex model\\ on
  planar maps}

\DeclareMathOperator{\tr}{tr}
\newcommand\ZZ{\mathbb{Z}}
\newcommand\CC{\mathbb{C}}

\author{Mireille Bousquet-M\'elou, Andrew Elvey Price and Paul Zinn-Justin}
\thanks{Paul Zinn-Justin was supported by ARC grant FT150100232.}
\address{LaBRI, CNRS, Universit\'e de Bordeaux, France \\
School of Mathematics \& Statistics, The University of Melbourne, Victoria 3010, Australia}

\newcommand{\Cgf}{{\sf C}}
\newcommand{\Dgf}{{\sf D}}
\newcommand{\Ggf}{{\sf G}}
\newcommand{\Hgf}{{\sf H}}
\newcommand{\Pgf}{{\sf P}}
\newcommand{\Qgf}{{\sf Q}}

\newcommand{\Rgf}{{\sf R}}

\newcommand{\Wgf}{{\sf W}}
\newcommand{\zs}{\mathbb{Z}}
\long\def\remint#1#2{%
\tikz[baseline=-4pt]%
{\node[outer sep=0pt,draw=black,fill=cyan!40!green!50!white,rectangle,rounded corners,align=left,text width=#1]{#2}}%
}
\newbox{\rembox}
\long\def\rem#1{\noindent\nobreak\hfil\penalty1000\hfilneg
\sbox{\rembox}{#1}%
\ifdim\wd\rembox>\textwidth\remint{\textwidth}{#1}\else\remint{}{#1}\fi%
}
\renewcommand{\d}{\mathrm{d}}

\newcommand\W{W^{(0)}}

\def\emm#1,{{\em #1}}
\newcommand{\gf}{generating function}

\newtheorem*{lem*}{Lemma}
\newtheorem{Theorem}{Theorem}

\newtheorem{Proposition}[Theorem]{Proposition}
\newtheorem{Definition}[Theorem]{Definition}

\newcommand{\beq}{\begin{equation}}
\newcommand{\eeq}{\end{equation}}

\renewcommand\th{\vartheta}
\usepackage[english]{babel}
\usepackage[utf8]{inputenc}

\begin{document}
\maketitle

\begin{abstract}
We address the enumeration of
  planar 4-valent maps equipped with an Eulerian orientation by two different methods, and
  compare the solutions we thus obtain. With the first method we enumerate these orientations as well as a restricted class which we show to be in bijection with general Eulerian orientations. The second method, based on the work of Kostov, allows us to enumerate these 4-valent orientations with a weight on some vertices, corresponding to the six vertex model. We prove that this result generalises both results obtained using the first method, although the equivalence is not immediately clear. 
\end{abstract}

\section{Introduction}
\begin{figure}[bht]
\centering
   \includegraphics[scale=0.7]{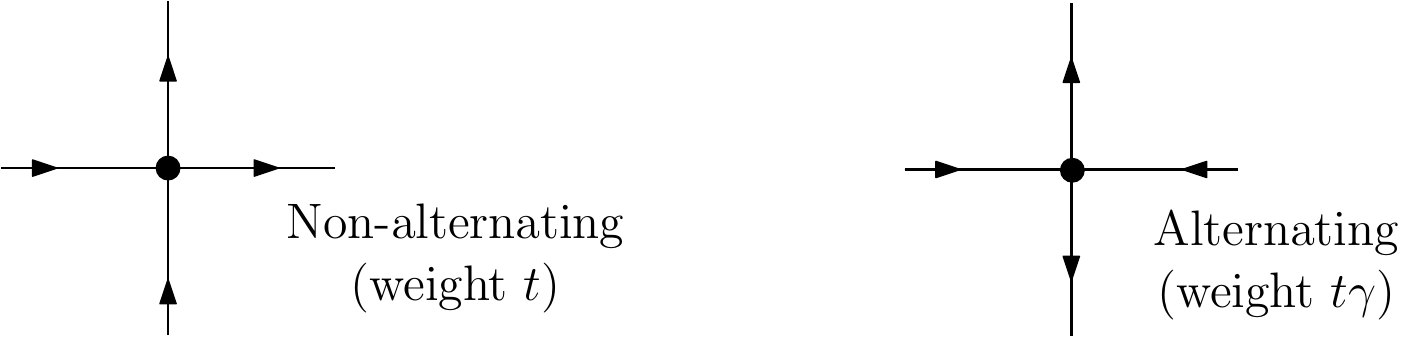} 
   \caption{The two types of vertices in the six-vertex
     model.}
   \label{fig:two vertex types}
\end{figure}

In 2000, Zinn-Justin \cite{artic10} and Kostov \cite{Kostov-6v}
studied the \emm six-vertex model on a random lattice,. In
combinatorial terms, this means counting rooted 4-valent (or: \emm
quartic,) planar maps equipped
with an \emm Eulerian orientation, of the edges: that is, every vertex
has equal in- and out-degree  (Figs.~\ref{fig:two vertex types} and~\ref{fig:M-AB}). Every vertex is
weighted $t$, and every \emm alternating, vertex gets an additional
weight $\gamma$, yielding a \gf\ $\Qgf(t,\gamma)$:
$$
Q(t, \gamma)=  \left( 2\,\gamma+2 \right) t+ \left( 9\,{\gamma}^{2}+16\,\gamma+10
 \right) {t}^{2}+ \left( 54\,{\gamma}^{3}+132\,{\gamma}^{2}+150\,
\gamma+66 \right) {t}^{3}+O( {t}^{4} ) 
.
$$
For instance, the 4 orientations accounting for the coefficient of $t$
are the following ones (the edge carrying a double arrow is the root edge, oriented canonically):
\begin{center}
  \includegraphics[scale=0.9]{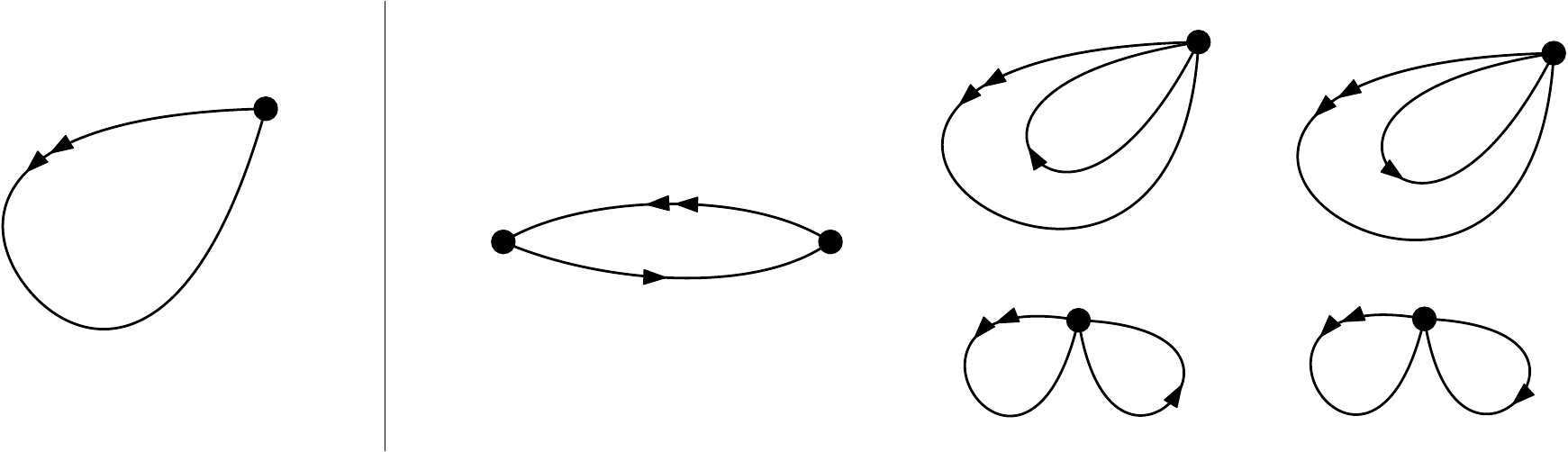}
\end{center}
Kostov solved this problem exactly, but the
form of his solution is quite complicated and its derivation is not entirely rigourous.  

Recently Bonichon, Bouquet-M\'elou, Dorbec and Pennarun
\cite{BoBoDoPe} posed the analogous 
problem of enumerating 
Eulerian orientations of general planar maps by edges
(where the vertex degree is not restricted).
They were followed by
Elvey Price and Guttmann who wrote an
intricate system of functional equations defining the associated
generating function \cite{elvey-guttmann17}. This allowed them to compute the number $g_n$ of Eulerian orientations
with $n$ edges for large values of $n$, and led them to a conjecture
on the asymptotic behaviour of $g_n$. In a similar way they conjectured the
asymptotic behaviour of the coefficients $q_n$ of $\Qgf(t,1)$,
counting Eulerian orientations of quartic maps, though this prediction
had already been in the physics
papers~\cite{Kostov-6v,artic10}. The conjectured asymptotic forms of
the sequences $(q_{n})_{n\geq0}$ and $(g_{n})_{n\geq0}$ led us to conjecture exact forms of the two sequences. These are the conjectures that we prove in the following two theorems.

\begin{Theorem}\label{thm:gen}
  Let $\Rgf(t)\equiv \Rgf$ be the unique formal power series with constant
  term $0$ satisfying
\beq \label{R0}
t= \sum_{n \ge 0} \frac 1 {n+1} {2n \choose n}^2\Rgf^{n+1}.
\eeq
Then the generating function of rooted planar Eulerian orientations, counted by
edges, is
\[
\Ggf(t)=\frac{1}{2}\Qgf(t,0)= \frac 1{4t^2}\left(
  t-2t^2-\Rgf(t)\right)= t+5t^2+33t^3+ \cdots.
\]
\end{Theorem}

\begin{Theorem}\label{thm:4}
  Let $\Rgf(t)\equiv \Rgf$ be the unique formal power series with constant
  term $0$ satisfying
\[
t= \sum_{n \ge 0} \frac{1}{n+1} {2n \choose n}{3n \choose n} \Rgf^{n+1}.
\]
Then the generating function of quartic rooted planar Eulerian orientations, counted by
vertices, is
\[
\Qgf(t,1)= \frac 1{3t^2}\left( t-3t^2-\Rgf(t)\right) =4t+35t^2+402 t^3+ \cdots.
\]
\end{Theorem}

The first step in our proof of Theorem \ref{thm:gen} is a 
bijection that relates general Eulerian orientations to quartic ones
\emm having no alternating vertex, (Section~\ref{sec:bij}). It implies that
$\Ggf(t)=\frac{1}{2}\Qgf(t,0)$. We then characterise the generating
functions $\Qgf(t,0)$ and $\Qgf(t,1)$ by a system of functional
equations using some new decompositions of planar maps. We then solve
these 
equations exactly (Section~\ref{sec:orientations}).
Details on this approach as well as basic definitions on planar maps can be found in~\cite{BM-EP18}.

In Section~\ref{sec:six} we use a different method to analyse the generating
function $\Qgf(t,\gamma)$ for general $\gamma$, following Kostov's
solution to the six-vertex model. We re-derive the first part of his
study, which yields a system of functional equations characterising
$\Qgf(t,\gamma)$, using a combinatorial argument.
We then follow Kostov's solution to these equations and fix a mistake,
which gives us a parametric expression  of $\Qgf(t,\gamma)$ in terms of the Jacobi theta function
\[
\th(z,q)=2\sin(z) q^{1/8} \prod_{n=1}^{\infty} (1-2\cos(2z)q^n+q^{2n})(1-q^n).
\]

\begin{Theorem}\label{thm:allgamma}
Write $\gamma=-2\cos(2\alpha)$, and   let $q(t,\gamma)\equiv q= t+ \left( 
6\,\gamma+6 \right) {t}^{2}+\cdots $ be the
unique formal power series in $t$ with constant
  term $0$ satisfying
\[
t=  \frac{\cos\alpha}{64\sin^3\alpha}
\left(
-\frac{\th(\alpha,q)\th'''(\alpha,q)}{\th'(\alpha,q)^2}+\frac{\th''(\alpha,q)}{\th'(\alpha,q)}
\right),
\]
 where all derivatives are with respect to the first variable. Moreover, define the series $\Rgf(t,\gamma)$ by
\[
\Rgf(t,\gamma)=\frac{\cos^2\alpha}{96\sin^4\alpha}
\frac{\th(\alpha,q)^2}{\th'(\alpha,q)^2}
\left(-\frac{\th'''(\alpha,q)}{\th'(\alpha,q)}
+\frac{\th'''(0,q)}{\th'(0,q)}\right).
\]
 Then the generating function of quartic rooted planar Eulerian orientations, counted by
vertices, with a weight $\gamma$ per alternating vertex is
\[
\Qgf(t,\gamma)= \frac{1}{(\gamma+2)t^2}\left( t-(\gamma+2)t^2-\Rgf(t,\gamma)\right).
\]
\end{Theorem}
It is not clear why, in  the two special cases
$\gamma=0,1$, this theorem is equivalent to Theorems \ref{thm:gen} and
\ref{thm:4} respectively. We prove this equivalence in
Section~\ref{sec:final}.
We conclude with a discussion on
further projects.


\begin{figure}[htb]
  \centering
   \scalebox{0.9}{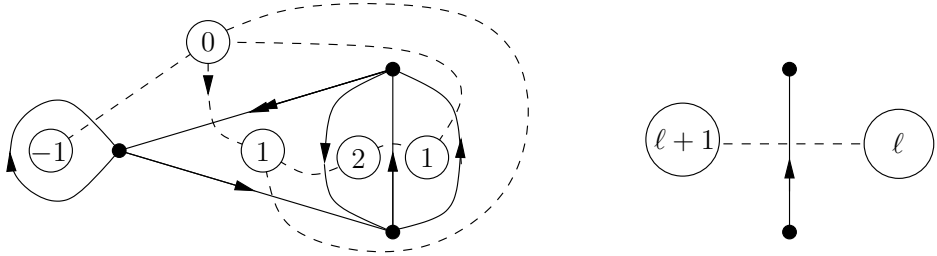}
  \caption{A rooted Eulerian orientation (solid edges;  the root edge
    is shown with a
    double arrow, and its orientation is chosen canonically) and the corresponding
    dual labelled map (dashed edges). The labelling rule is shown on
    the right.}
  \label{fig:duality}
\end{figure}

\section{Bijection for general Eulerian orientations}
\label{sec:bij}
The first step in  the enumeration of Eulerian orientations is a
simple bijection, introduced in~\cite{elvey-guttmann17}, to certain
{\em labelled maps}. 
\begin{Definition}\label{def:labelled-map}
  A {\em labelled map}\/ is a rooted planar map with integer labels on
  its vertices,  such that  adjacent labels differ by $1$ and
  the root edge is labelled from $0$ to $1$.
\end{Definition}

The bijection  is illustrated in
Figure~\ref{fig:duality}. The idea is that an Eulerian orientation of
edges of a map determines a height function on the vertices of its
dual. A restriction of this bijection shows that {\em quartic}\/
Eulerian orientations are in bijection with labelled {\em
  quadrangulations} (every face has degree 4). 

\begin{figure}[b!]
  \centering
  \includegraphics[scale=0.7]{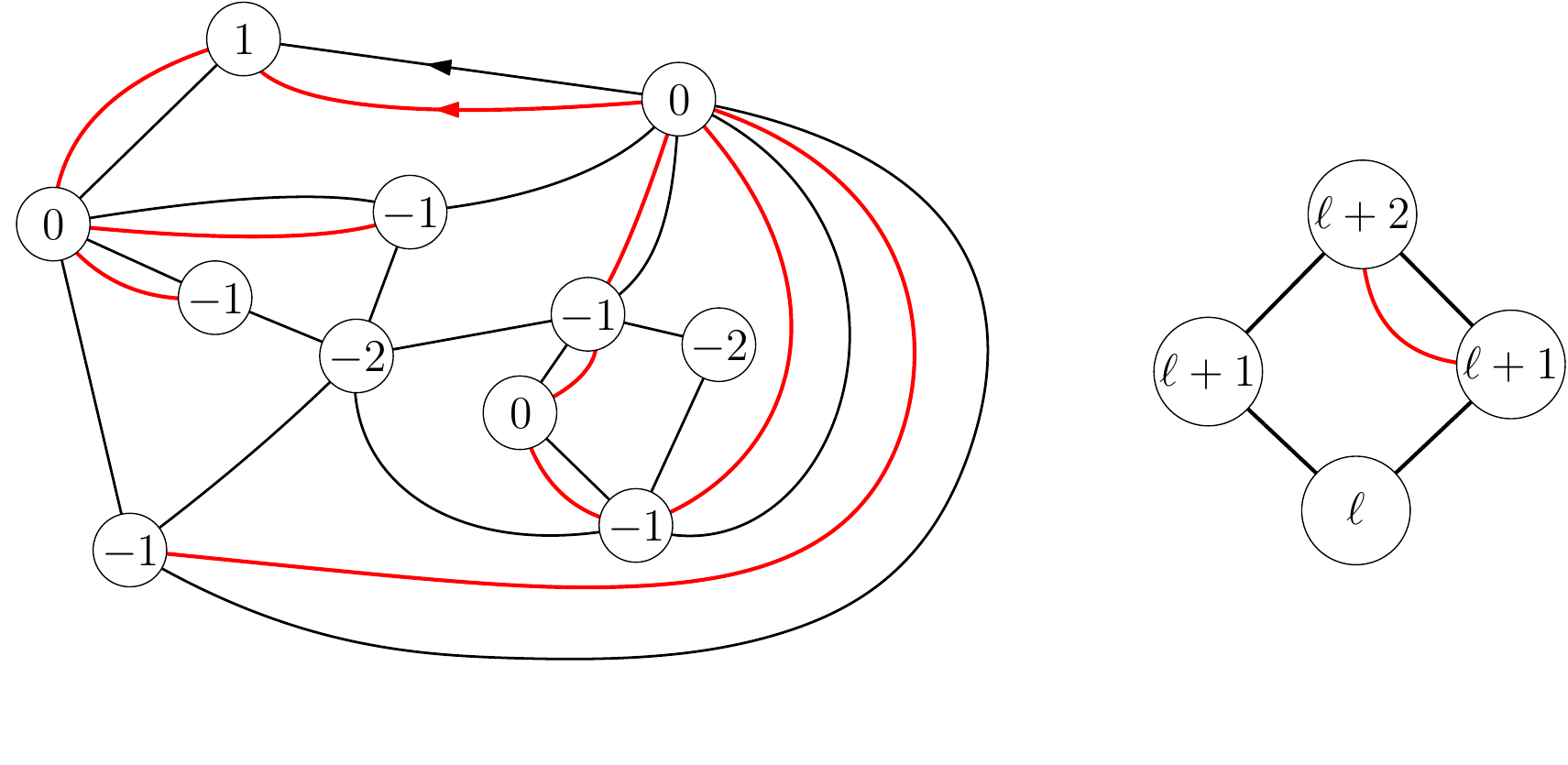}
  \caption{A  labelled quadrangulation $Q$ (black edges) and the corresponding
    labelled map $L$ (red edges). The rule for drawing red edges is shown on
    the right. Note that the two local minima of $Q$, both labelled
    $-2$, disappear in the construction.}
  \label{fig:M-AB}
\end{figure}

For the next step, we use a bijection of Miermont, Ambj{\o}rn and Budd
\cite{miermont2009tessellations,ambjorn2013trees} to show that
labelled maps having $n$ edges are in 1-to-2 correspondence  with \emm colourful,
labelled quadrangulations having $n$ faces. By \emm colourful, we
mean that  each face has three distinct labels, or equivalently, that
the corresponding quartic orientation has no alternating vertex.
This bijection generalizes a bijection
of~\cite{chassaing-schaeffer}, and is
illustrated in Figure \ref{fig:M-AB}. It implies that $\Ggf(t)=\frac{1}{2}\Qgf(t,0)$.

In fact,  this pair of bijections allows us to understand the more
general series $\Qgf(t,\gamma)$ as a generalisation of $\Ggf(t)$ in
terms of   Eulerian  \emm partial, orientations. These are planar maps
in which \emm some, edges are oriented, in such a way  that each
vertex has equal in- and out-degree. 

\begin{Proposition}
  The series $\Qgf(t,\gamma)$ counting rooted  quartic Eulerian orientations also counts
rooted   Eulerian partial orientations with a weight $t$ per edge and
an additional weight $\gamma$ per undirected edge (the root edge may be
undirected or directed in either direction).
\end{Proposition}
Here are the 4 partial orientations that account for the coefficient
$2+2\gamma$ of $t$ in $\Qgf(t,\gamma)$:
\begin{center}
  \includegraphics[scale=0.9]{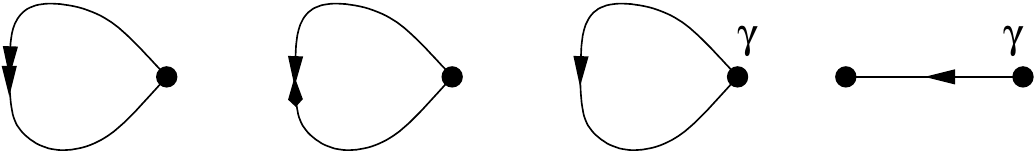}
\end{center}

\section{The number of planar Eulerian orientations}
\label{sec:orientations}
In this section we give a very brief summary of our solutions to the cases $\Qgf(t,0)$ and $\Qgf(t,1)$. The full details are in \cite{BM-EP18}.
We define three classes $\mathcal P, \mathcal D$ and $\mathcal C$ of
labelled quadrangulations and decompose them recursively, using in
particular a new contraction operation.  We thus
characterise the series $\Qgf(t,0)=2\Ggf(t)$ as follows.

\begin{Proposition}\label{thm:systemG}
  There is a unique $3$-tuple of series, denoted $\Pgf(t,y)$, $\Cgf(t,x,y)$
and $\Dgf(t,x,y)$, belonging respectively to $\mathbb{Q}[[y,t]]$,
  $\mathbb{Q}[x][[y,t]]$ and $\mathbb{Q}[[x,y,t]]$, and satisfying the
  following equations:
\begin{align*}
\Pgf(t,y)&=\frac{1}{y}[x^1]\Cgf(t,x,y),\\
\Dgf(t,x,y)&=\frac{1}{1-\Cgf\left(t,\frac{1}{1-x},y\right)},\\
\Cgf(t,x,y)&=xy[x^{\geq0}]\left(\Pgf(t,tx)\Dgf\left(t,\frac{1}{x},y\right)\right),
\end{align*}
together with the initial condition $\Pgf(t,0)=1$ (the operator
$[x^{\geq0}]$ extracts all monomials in which the exponent of $x$ is non-negative).

The generating function $\Qgf(t,0)$ is related to these series by the equation
\[
\Qgf(t,0)=[y^1]\Pgf(t,y)-1.
\]
\end{Proposition}
We then solve this system of equations as follows, writing $\Pgf$, $\Cgf$ and $\Dgf$ in terms of the series $\Rgf(t)$ defined in Theorem \ref{thm:gen}:
\[
t\Pgf(t,ty)=\sum_{n\ge 0}\sum_{j=0}^n \frac 1{n+1} {2n\choose n}{2n-j
  \choose n} y^j \Rgf^{n+1},
\]
\[
\Cgf(t,x,ty) =1 -\exp\left( 
- \sum_{n\ge 0} \sum_{j=0}^n \sum_{i=0}^{n}
  \frac 1 {n+1} {2n-i \choose n} {2n-{j}\choose n}  x^{i+1}  y^{j+1}\Rgf^{n+1} \right),
\]
\[
\Dgf(t,x,ty) =\exp\left( \sum_{n\ge 0} \sum_{j=0}^n \sum_{i\ge0}
  \frac 1 {n+1} {2n-j \choose n} {2n+i+1\choose n}  x^i   y^{j+1}\Rgf^{n+1}\right).
\]
Theorem \ref{thm:gen} then follows from the equation
$2\Ggf(t)=\Qgf(t,0)=[y^1]\Pgf(t,y)-1.$  We first obtained the solution
 using a guess and check approach, but we now have a constructive way
 of deriving it from Proposition~\ref{thm:systemG}.

We have a similar proof of Theorem \ref{thm:4}: we characterise
$\Qgf(t,1)$ using a  system of functional equations, which we then solve.

\section{The  matrix integral approach to the six-vertex model}
\label{sec:six}
Following \cite{artic10} and \cite{Kostov-6v}, we introduce the following
matrix integral:
\begin{equation}\label{eq:defZ}
Z_N=\int \d X \d X^\dagger \exp\left[N\tr\left(
-XX^\dagger + t X^2X^\dagger{}^2+{\gamma t\over 2}(XX^\dagger)^2\right)\right]
\end{equation}
where integration is over $N\times N$ complex matrices,
and $X^\dagger$ denotes the conjugate transpose of $X$. Then
\begin{equation}\label{eq:F}
2t\frac{\partial}{\partial t}\log Z_N = \sum_{g\ge 0} \Qgf^{(g)}(t,\gamma) N^{2-2g},
\end{equation}
where each series $\Qgf^{(g)}(t,\gamma)$ is the genus $g$ analogue of $\Qgf(t,\gamma)\equiv\Qgf^{(0)}(t,\gamma)$.

Extracting the series $\Qgf(t,\gamma)$ from Kostov's solution \cite{Kostov-6v} by using \eqref{eq:F} directly is not easy, so instead we start with a combinatorial interpretation of Kostov's work in which $\Qgf(t,\gamma)$ appears naturally.

We first convert the matrix integral~\eqref{eq:defZ} into another
integral, this time involving three matrices, which can be understood
in terms of cubic maps. Then, a standard first step is
to derive from such integrals ``loop equations'' relating certain \emm correlation
functions,. In fact we also have  direct  combinatorial
proofs of these equations in terms of certain families of partially
oriented  maps.

\begin{Proposition}\label{6vertex_combi} There is a unique pair of
  series, denoted $\Wgf(t,\omega,x)\equiv\Wgf(x)$ and
  $\Hgf(t,\omega,x,y)\equiv \Hgf(x,y)$, belonging respectively to
  $\mathbb{Q}(\omega)[x][[t]]$ and $\mathbb{Q}(\omega)[x,y][[t]]$ and satisfying the equations
\begin{align*}\Wgf(x)&=x^2t\Wgf(x)^{2}+\omega xt\Hgf(0,x)+\omega^{-1}xt\Hgf(x,0)+1\\
\Hgf(x,y)&=\Wgf(x)\Wgf(y)+\frac{\omega}{y}\left(\Hgf(x,y)-\Hgf(x,0)\right)+\frac{\omega^{-1}}{x}\left(\Hgf(x,y)-\Hgf(0,y)\right).\end{align*}
The series $\Qgf(t,\gamma)$ is given by 
\[
  \Qgf\left(t,\omega^{2}+\omega^{-2}\right)=\Hgf(t,\omega,0,0)-1=\frac{1}{t(\omega+\omega^{-1})}[x^1]\Wgf(x)
  -1.
\]
\end{Proposition}

\begin{figure}[ht]
\centering
   \includegraphics[scale=1]{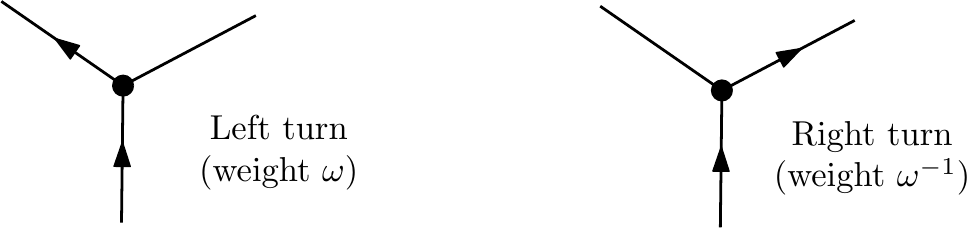} 
   \caption{The two vertex types allowed as non-root vertices in the
     Eulerian partial orientations counted by $\Wgf(t,\omega,x)$ and
     $\Hgf(t,\omega,x,y)$.}
   \label{fig:two_A_vertex_types}
\end{figure}

Each of the two series $\Wgf$ and $\Hgf$ counts some class of
rooted  Eulerian  \emm partial, orientations in which each non-root
vertex is one of the two types shown in Figure
\ref{fig:two_A_vertex_types}, with weights $\omega$ and $\omega^{-1}$
as shown, and $t$ counts 
undirected edges. The series $\Wgf$ and $\Hgf$ differ in the weight and allowed type of the root vertex. The equations relating $\Hgf$ and $\Wgf$ follow from contracting the root edge, while the relation to $\Qgf$ follows from contracting all of the undirected edges.

%
To solve these equations we convert the series back to Kostov's setting via the transformation
\[\W(x)=\frac{1}{x}\,\Wgf\!\left(\frac{1}{x}\right),\]
and similar transformations for $\Hgf(x,y)$, $\Hgf(x,0)$ and
$\Hgf(0,y)$. We reinterpret these transformed series as complex
analytic functions of $x$, with $t$ a fixed small real number and
$\omega$ fixed, with $|\omega|=1$. In order to solve these equations one first proves a
technical lemma (the ``one-cut lemma'') which states that
the
function $\W(x)$ is analytic in $x$ except on a single cut $[x_1,
x_2]$ on the positive real line. After some
algebraic manipulation of the equations in Proposition 
\ref{6vertex_combi}, 
we arrive at the equation
\[0=\W(x+i0)+\W(x-i0)-\frac{x}{t}+\omega^{-2}\W(\omega^{-1}-\omega^{-2}x)+\omega^2\W(\omega-\omega^2x),\]
for $x\in\mathbb{R}$ on the cut of $\W(x)$ (see~\cite[Eq.~(3.19)]{Kostov-6v}).  As Kostov explains, the function
\begin{align*}
U(x):=x\omega \W\left(\frac{1}{\omega+\omega^{-1}}+i\omega x\right)&+x\omega^{-1}\W\left(\frac{1}{\omega+\omega^{-1}}-i\omega^{-1}x\right)\\
&+\frac{ix^2}{t(\omega^2-\omega^{-2})}  -\frac{x}{t(\omega+\omega^{-1})^2}
\end{align*}
is uniquely defined by the fact that $U(x)$ is holomorphic in
$\mathbb{C}$ minus the two cuts $( i\omega)^{\pm 1}[x'_1,x'_2]$, where
$x_i'$ is a translate of $x_i$ by an explicit real constant, along with the following three equations
\begin{equation}\label{eq:glue}
U(i\omega(x\pm i0))=U(-i\omega^{-1}(x\mp i0)),
\qquad
x\in (x'_1,x'_2)
\end{equation}
\begin{equation}\label{eq:exp}
U(x)=\frac{i}{t(\omega^2-\omega^{-2})} x^2 - \frac{1}{t(\omega+\omega^{-1})^2}x
+O(1/x) \qquad \hbox{as }
x\rightarrow \infty,
\end{equation}
\begin{equation}\label{eq:norm}
\oint_{\mathcal C} \frac{dx}{2\pi x} U(x)=1 ,
\end{equation}
where $\mathcal C$ surrounds the cut $(i\omega)^{-1}[x'_1,x'_2]$ anticlockwise.
 
Note that by expanding $U(x)$ at infinity further than \eqref{eq:exp}, i.e., $U(x)=\sum_{i=-2}^\infty U_i x^{-i}$,
we can extract from $U(x)$ the same information as from $\W(x)$. In particular,
\begin{equation}\label{eq:U1}
U_1=
1-(\omega+\omega^{-1})[x^{-2}]\W(x)=1-t(\omega+\omega^{-1})^2\left(1+\Qgf\!\left(t,\omega^{2}+\omega^{-2}\right)\right).
\end{equation}

\subsection*{Solution in terms of theta functions}
We now provide a parametric expression for $U(x)$, following \cite{Kostov-6v}.
We will first parametrise the domain of analyticity of $U(x)$. Let
$\th$ denote the classical Jacobi theta function $\theta_1$:
\begin{align}\label{theta-def}
\th(z)\equiv\th(z,q):=\theta_1(z;\tau)&= 2\sin(z) q^{1/8} \prod_{n=1}^{\infty} (1-2\cos(2z)q^n+q^{2n})(1-q^n),\\
&= -i \sum_{n\in \zs} (-1)^n e^{(n+1/2)^2\pi i\tau+(2n+1)iz},
\end{align}
where $q=e^{2\pi i \tau}$ and $\tau$ has positive imaginary
  part.  Define the mapping $x:\mathbb{C}\to\mathbb{C}$ by
\[
x(z) = x_0 \frac{\th(z+\alpha)}{\th(z)},
\]
where $\alpha$ is chosen so that $\omega=ie^{-i\alpha}$, which gives $\gamma=\omega^{2}+\omega^{-2}=-2\cos(2\alpha)$. The quantities $x_0$ and $\tau$ will
be determined later. Note that $x(z)$ is a meromorphic function whose poles form the lattice $\pi\tau \ZZ+\pi\ZZ$.

From the identities
  \[
    \th(z+\pi)=-\th(z) \quad \hbox{and} \quad \th(z+\pi\tau)=-e^{-i\pi
      \tau -2i z}\th(z),
  \]
  we obtain the pseudo-periodic identities:
\[
x(z+\pi)=x(z),
\qquad \hbox{and} \qquad 
x(z+\pi\tau)=
e^{-2i\alpha}
x(z) .
\]
The former identity implies that $x$ is a meromorphic function on the cylinder $C=\mathbb{C}/(\pi\mathbb{Z})$. 
Then, as explained in \cite{Kostov-6v}, property~\eqref{eq:glue} will be satisfied for certain $x_{1}',x_{2}'\in\mathbb{R}$ provided that the following holds: There is some meromorphic function $V(z)$ on the complex torus $T=C/(\pi\tau \ZZ)=\CC/(\pi\ZZ+\pi\tau \ZZ)$ and some fundamental domain $\hat{T}\subset C$ of $T$ containing $0$ such that $U(x(z))=V(z)$ for $z\in\hat{T}$. The restriction of $x$ to $\hat{T}$ sends each of the two boundaries of $\hat{T}$ to one of the two cuts in the domain of analyticity of $U(x)$.

Furthermore, because of the analyticity properties of $U$, 
the only singularity of $V(z)=U(x(z))$ comes from the double pole of $U(x)$
at $x=\infty$, i.e., $z=0$. Since $V(z)$ is meromorphic on $T$ and its only singularity is a double pole at $z=0$, it must be a linear transformation of the Weierstrass function:
\[
U(x(z))=A+B\wp(z),
\qquad
\wp(z)=\frac{1}{z^2}+\sum_{(m,n)\ne (0,0)} \left(\frac{1}{(z+\pi(m+n\tau))^2}-\frac{1}{\pi^2(m+n\tau)^2}\right).
\]
The parameters $\tau,x_0,A,B$ are determined by the expansion of $U$ at infinity \eqref{eq:exp} and
the normalization condition~\eqref{eq:norm}.

The three terms of expansion~\eqref{eq:exp} provide three equations which determine $x_{0}$, $B$ and $A$ in terms of $\alpha$ and $\tau$. We ignore the equation coming from the constant term, since this only determines $A$, which plays no role in any further calculations. We are left with the equations:
\begin{equation*}
B=\frac{\cos\alpha}{\sin^2\alpha}
\th_3^4(0)
\frac{\th^2(\alpha)}{\th'^2(\alpha)}
,\qquad
x_0=\frac{\cos\alpha}{2\sin\alpha}\frac{\th'(0)}{\th'(\alpha)}
\end{equation*}
where $\th_3(0)=\sum_{n=-\infty}^{+\infty} q^{n^2/2}$.
The integral \eqref{eq:norm} can be computed; fixing a mistake in~\cite[App.~B.2]{Kostov-6v} results in
a massive simplification:
\begin{equation}\label{eq:b}
t=
  \frac{\cos\alpha}{64\sin^3\alpha}
\left(
-\frac{\th(\alpha)\th'''(\alpha)}{\th'(\alpha)^2}+\frac{\th''(\alpha)}{\th'(\alpha)}
\right).
\end{equation}
The last equation should be understood as an implicit equation for $q=e^{2\pi i\tau}$ as a function of
$t$; if we want to return to formal power series, 
then it determines $q$ uniquely once we require $t\sim q$ around $0$,
as claimed in Theorem~\ref{thm:allgamma}.

Finally, by expanding $U(x)$ one order further, one finds
\begin{equation}\label{eq:W1}
t U_1=
\frac{\cos^2\alpha}{96\sin^4\alpha}
\frac{\th(\alpha)^2}{\th'(\alpha)^2}
\left(-\frac{\th'''(\alpha)}{\th'(\alpha)}
+\frac{\th'''(0)}{\th'(0)}\right).
\end{equation}
Theorem \ref{thm:allgamma} then follows due to \eqref{eq:U1}, writing $\Rgf(t,\gamma)=t U_1$, where $\gamma=-2\cos(2\alpha)$.

\section{Relationships between results and further problems}
\label{sec:final}

\subsection*{The cases $\gamma=0$ and $\gamma=1$}
We now describe how our formulas for $\Qgf(t,0)$ and $\Qgf(t,1)$ in
Theorems \ref{thm:gen} and \ref{thm:4} can be derived from our general
formula for $\Qgf(t,\gamma)$ in Theorem \ref{thm:allgamma}. It
suffices to show that the series $\Rgf(t,\gamma)$ coincides with the
series $\Rgf(t)$ in each case. In sight of~\eqref{R0}, for $\gamma=0$ this is equivalent to
\beq\label{toprove0}
  t=\Rgf(t,0)\ {}_2F_1\left(1/2,1/2;2|16\, \Rgf(t,0)\right),
\eeq
where we use the standard hypergeometric notation, 
while for $\gamma=1$, it is equivalent to
\beq\label{toprove1}
t=\Rgf(t,1)\ {}_2F_1(1/3,2/3;2|27\, \Rgf(t,1)).
\eeq
In each case, we 
rewrite both sides of the equation as series in $q$ using the
expressions in Theorem \ref{thm:allgamma}, noting that $\gamma=0$ corresponds to $\alpha=\pi/4$ and $\gamma=1$ corresponds to $\alpha=\pi/3$. We then introduce the following parametrisation of $q$,
as first suggested in a slightly different language by Ramanujan \cite{Berndt-Rama}:
\begin{gather*}
  q
  =\exp\left(-\frac{\pi}{\sin(\pi a)}
\frac{\tilde A(w)}{A(w)}\right),
\\
A(w)={}_2F_1(a,1-a;1|w),
\qquad
\tilde A(w)=A(1-w),
\end{gather*}
where $a$ is a constant to be specified shortly, and $w$ is the new parameter.
The hypergeometric series $A$ and $\tilde A$ 
satisfy the {\em same}\/   differential equation:
\begin{equation}\label{eq:hyper}
w(1-w)\frac{d^2A}{dw^2}+(1-2w)\frac{dA}{dw}-a(1-a)A=0.
\end{equation}
Writing $q=e^{2\pi i\tau}$,  it is not hard to prove that $\tau(w)$
satisfies the following equation:
\begin{equation}\label{eq:tauid}
2\pi i\, w(1-w)\frac{d\tau}{dw}\,
A(w)^2=1.
\end{equation}
Finally, the usual differentiation formula for hypergeometric series yields
\begin{equation}\label{eq:derivid}
(1-w)\frac{dA}{dw}
=a(1-a){}_2F_1(a,1-a;2|w).
\end{equation}

The functional inverse $w(\tau)$ is known explicitly in the four cases
$a=1/2,1/3,1/4,1/6$ \cite{borwein1991cubic}; the
first two will be relevant to us:
\[
w(\tau)=
\left\{
\begin{aligned}
&\left(\frac{C}{A}\right)^2,& A&=\sum_{m,n\in\ZZ}q^{m^2+n^2},& C&=\sum_{m,n\in\ZZ+1/2}q^{m^2+n^2},
&
a&=1/2
\\
&\left(\frac{C}{A}\right)^3,& A&=\sum_{m,n\in\ZZ}q^{m^2+mn+n^2},& C&=\sum_{m,n\in\ZZ+1/3} q^{m^2+mn+n^2},
&
a&=1/3
\end{aligned}
\right.
\]
where in addition, $A$ coincides with  $A(w(\tau))$. In both cases, one has the following identity, which follows from the product
definition \eqref{theta-def} of $\th$:
\begin{equation}\label{eq:Ath}
A=\tan\alpha \, \frac{\th'(\alpha,q)}{\th(\alpha,q)},
\qquad
\left\{
\begin{aligned}
  a&=1/2,& \alpha&=\pi/4,
  \\
  a&=1/3,& \alpha&=\pi/3.
\end{aligned}
\right.
\end{equation}
These series are  entries 
\href{https://oeis.org/A004018}{A004018}
and
\href{https://oeis.org/A004016}{A004016} in the OEIS~\cite{oeis},
respectively, and the identities above can be found there. We will also use the heat equation
\beq\label{heat}
  \th''(z,q):=\frac{\partial ^2 \th}{\partial z^2} (z,q) =-\frac{4}{i
    \pi }\frac{d \th}{d\tau}(z,q),
  \eeq
  which allows us to express all $\tau$-derivatives of $\th$ in terms
  of $z$-derivatives.

We compute $\Rgf=t U_1$ from \eqref{eq:W1} using the function $\eta(x)=\prod_{n\ge1}(1-x^n)$.
In the case $\alpha=\pi/4$ we have the equations $-\th'''(0)/\th'(0)=1+24q \frac{\eta'(q)}{\eta(q)}$
(\href{https://oeis.org/A006352}{A006352})
and $-\th'''(\pi/4)/\th'(\pi/4)=1-24q\frac{\eta'(-q)}{\eta(-q)}$ (\href{https://oeis.org/A143337}{A143337}), so
\[48\Rgf A^2=\frac{\th'''(0)}{\th'(0)}-\frac{\th'''\left(\frac{\pi}{4}\right)}{\th'\left(\frac{\pi}{4}\right)}=-24q\frac{\eta'(q)}{\eta(q)}-24q\frac{\eta'(-q)}{\eta(-q)}=3\sum_{n\ge0}\frac{(2n+1)q^{2n}}{1-q^{4n+2}}
=3C^2,\]
and therefore $\Rgf=\frac{1}{16}w$. The last equality above is from \href{https://oeis.org/A008438}{A008438}. For $\alpha=\pi/3$, we show that $\Rgf=\frac{1}{27}w$ in a similar way, although the proof is more complicated.

Using the parametrisation in terms of $w$, we are now in a position to
prove the identities~\eqref{toprove0} and~\eqref{toprove1}, with $t$ and
$\Rgf(t,\gamma)$ given by Theorem~\ref{thm:allgamma}.
\begin{align*}
t&=\frac{1}{8\sin^2\alpha}\frac{1}{2\pi i}\, A^{-2}\frac{dA}{d\tau} &&\text{from \eqref{eq:b}, \eqref{eq:Ath}, \eqref{heat}}
\\
&=\frac{1}{8\sin^2\alpha}\, w(1-w)\frac{dA}{dw}&&\text{from \eqref{eq:tauid}}
\\
&=\frac{a(1-a)}{8\sin^2\alpha}\,w\  {}_2F_1(a,1-a;2|w)&&\text{from \eqref{eq:derivid}}.
\end{align*}
The desired results then follow as $\Rgf=\frac{a(1-a)}{8\sin^2\alpha}w$ in both cases.


\subsection*{Generalisations and further questions}

When $\gamma=0$ or $\gamma=1$, we have obtained two different
parametric expressions of $\Qgf(t,\gamma)$: one in terms of
hypergeometric series, the other in terms of theta functions. Is there
an analogue of the hypergeometric version for general $\gamma$? We
have generalised the equations of Proposition~\ref{thm:systemG} to
include a weight  $\gamma$, but so far we have been unable to solve
them for general $\gamma$. 

In the case of general Eulerian orientations, we are interested in one
other natural generalisation of $\Ggf(t)$: the generating function
$\Ggf(t,z)$ which counts Eulerian orientations by edges ($t$) {\em
  and}\/ vertices ($z$).
Through the sequence of bijections of Section~\ref{sec:bij}, the
number of vertices in an Eulerian orientation is 
 the number of clockwise faces in the corresponding quartic
 orientation. 
 This is not a very natural quantity from the matrix integral
 perspective,
and it is not clear how to generalise the equations of Theorem
\ref{6vertex_combi} to include $z$. Nevertheless, we have generalised
the equations of Theorem \ref{thm:systemG} to include  $z$ (and in
fact $z$ and $\gamma$ simultaneously). In the specific cases
$\gamma=0$ and $\gamma=1$ we can solve these equations, thus generalising Theorems \ref{thm:gen} and \ref{thm:4}. 

\bibliographystyle{amsalphahyper}
\bibliography{biblio}

\end{document}

%% file: duality-small.pdf_t
\begin{picture}(0,0)%
\includegraphics{duality-small.pdf}%
\end{picture}%
\setlength{\unitlength}{4144sp}%
\begingroup\makeatletter\ifx\SetFigFont\undefined%
\gdef\SetFigFont#1#2#3#4#5{%
  \reset@font\fontsize{#1}{#2pt}%
  \fontfamily{#3}\fontseries{#4}\fontshape{#5}%
  \selectfont}%
\fi\endgroup%
\begin{picture}(6118,1669)(600,-3038)
\put(2251,-2423){\makebox(0,0)[lb]{\smash{{\SetFigFont{12}{14.4}{\familydefault}{\mddefault}{\updefault}{\color[rgb]{0,0,0}1}%
}}}}
\put(2880,-2469){\makebox(0,0)[lb]{\smash{{\SetFigFont{12}{14.4}{\familydefault}{\mddefault}{\updefault}{\color[rgb]{0,0,0}2}%
}}}}
\put(3331,-2468){\makebox(0,0)[lb]{\smash{{\SetFigFont{12}{14.4}{\familydefault}{\mddefault}{\updefault}{\color[rgb]{0,0,0}1}%
}}}}
\put(766,-2424){\makebox(0,0)[lb]{\smash{{\SetFigFont{12}{14.4}{\familydefault}{\mddefault}{\updefault}{\color[rgb]{0,0,0}$-1$}%
}}}}
\put(6436,-2371){\makebox(0,0)[lb]{\smash{{\SetFigFont{12}{14.4}{\familydefault}{\mddefault}{\updefault}{\color[rgb]{0,0,0}$\ell$}%
}}}}
\put(4871,-2343){\makebox(0,0)[lb]{\smash{{\SetFigFont{12}{14.4}{\familydefault}{\mddefault}{\updefault}{\color[rgb]{0,0,0}$\ell +1$}%
}}}}
\put(1891,-1696){\makebox(0,0)[lb]{\smash{{\SetFigFont{12}{14.4}{\familydefault}{\mddefault}{\updefault}{\color[rgb]{0,0,0}0}%
}}}}
\end{picture}%

%% file: sixvertexfinal.bbl
\def\cprime{$'$} \def\cprime{$'$}
\providecommand{\bysame}{\leavevmode\hbox to3em{\hrulefill}\thinspace}
\providecommand{\MR}{\relax\ifhmode\unskip\space\fi MR }
\providecommand{\MRhref}[2]{%
  \href{http://www.ams.org/mathscinet-getitem?mr=#1}{#2}
}
\providecommand{\href}[2]{#2}
\begin{thebibliography}{BBMDP17}

\bibitem[AB13]{ambjorn2013trees}
J.~Ambj{\o}rn and T.~G. Budd, \emph{Trees and spatial topology change in causal
  dynamical triangulations}, J. Phys. A \textbf{46} (2013), no.~31, 315201, 33.
  \MR{3090757}.

\bibitem[BB91]{borwein1991cubic}
J.~M. Borwein and P.~B. Borwein, \emph{A cubic counterpart of jacobi's identity
  and the agm}, Transactions of the American Mathematical Society (1991),
  691--701.

\bibitem[BBMDP17]{BoBoDoPe}
N.~Bonichon, M.~Bousquet-M\'elou, P.~Dorbec, and C.~Pennarun, \emph{On the
  number of planar {E}ulerian orientations}, European J. Combin. \textbf{65}
  (2017), 59--91, \href{https://arxiv.org/abs/1610.09837}{arXiv:1610.09837}.
  \MR{3679837}.

\bibitem[Ber98]{Berndt-Rama}
B.~C. Berndt, \emph{Ramanujan's notebooks --- part v}, Springer Publishing
  Company, Incorporated, 1998,
  {\scriptsize\href{http://dx.doi.org/10.1007/978-1-4612-1624-7}{\path{doi:10.1007/978-1-4612-1624-7}}}.

\bibitem[BMEP]{BM-EP18}
M.~Bousquet-M\'elou and A.~Elvey~Price, \emph{The generating function of planar
  eulerian orientations},
  \href{https://arxiv.org/abs/1803.08265}{arXiv:1803.08265}.

\bibitem[CS04]{chassaing-schaeffer}
P.~Chassaing and G.~Schaeffer, \emph{Random planar lattices and integrated
  super{B}rownian excursion}, Probab. Theory Related Fields \textbf{128}
  (2004), no.~2, 161--212,
  \href{https://arxiv.org/abs/math/0205226}{arXiv:math/0205226}. \MR{MR2031225
  (2004k:60016)}.

\bibitem[EPG18]{elvey-guttmann17}
A.~Elvey~Price and A.~J. Guttmann, \emph{Counting planar {E}ulerian
  orientations}, Europ. J. Combinatorics \textbf{71} (2018), 73–98,
  \href{https://arxiv.org/abs/1707.09120}{arXiv:1707.09120}.

\bibitem[Inc]{oeis}
OEIS~Foundation Inc., \emph{The on-line encyclopedia of integer sequences},
  \href{http://oeis.org}{http://oeis.org}.

\bibitem[Kos00]{Kostov-6v}
I.~Kostov, \emph{Exact solution of the six-vertex model on a random lattice},
  Nuclear Phys. B \textbf{575} (2000), no.~3, 513--534,
  \href{http://arxiv.org/abs/hep-th/9911023}{\path{arXiv:hep-th/9911023}}.
  \MR{MR1762323 (2001f:82019)}.

\bibitem[Mie09]{miermont2009tessellations}
G.~Miermont, \emph{Tessellations of random maps of arbitrary genus}, Ann. Sci.
  {\'E}c. Norm. Sup{\'e}r.(4) \textbf{42} (2009), no.~5, 725--781.

\bibitem[ZJ00]{artic10}
P.~Zinn-Justin, \emph{The six-vertex model on random lattices}, Europhys. Lett.
  \textbf{50} (2000), no.~1, 15--21,
  \href{http://arxiv.org/abs/cond-mat/9909250}{\path{arXiv:cond-mat/9909250}}.
  \MR{MR1747830 (2001a:82046)}.

\end{thebibliography}
